\documentclass[12pt]{article}

\usepackage{fleqn,amssymb,amsmath,amscd,epsfig,theorem} %
\newtheorem{theorem}{Theorem}

\newtheorem {proposition}{Proposition}
\newtheorem {corollary}{Corollary}
\theorembodyfont{\normalfont}
\newtheorem {definition}{Definition}
\newtheorem {example}{Example}

\newtheorem {remark}{Remark}
\newcommand{\beq}{\begin{equation}}
\newcommand{\eeq}{\end{equation}}
\newcommand{\Leq}[1]{\label{#1}\end{equation}}
\newcommand{\beqn}{\begin{eqnarray}}
\newcommand{\eeqn}{\end{eqnarray}}
\newcommand{\beqno}{\begin{eqnarray*}}
\newcommand{\eeqno}{\end{eqnarray*}}
\renewcommand {\l}{\left}
\newcommand {\ri}{\right}
\newcommand {\eh}{{\textstyle \frac{1}{2}}}
\newcommand {\ar}{\rightarrow}
\newcommand {\sign}{{\rm sign}}

\newcommand {\Id}{{\rm Id}}

\newcommand {\rank}{{\rm rank}}

\newcommand {\bN}{{\mathbb N}}
\newcommand {\bR}{{\mathbb R}}
\newcommand {\bZ}{{\mathbb Z}}
\newcommand{\rstr}{{\upharpoonright}}
\newcommand{\cF}{{\cal F}} 
\newcommand{\cH}{{\cal H}}
\newcommand{\cO}{{\cal O}} 
\newcommand{\cP}{{\cal P}}
\newcommand{\cS}{{\cal S}} 
\newcommand{\NN}{\nonumber}

\newcommand{\hH}{{\hat{H}\,}}
\newcommand{\hM}{{\hat{M}\,}}
\newcommand{\qmbox}[1]{\quad\mbox{#1}\quad}
\renewcommand {\max}{{{\rm max}}}

\newcommand{\ul}{{\underline{l}}}

\newcommand {\q}{{\vec{q}}}
\newcommand {\p}{{\vec{p}}}
\newcommand {\Q}{{\vec{Q}}}
\renewcommand {\P}{{\vec{P}}}
\newcommand {\pq}{{(\p,\q)}}
\newcommand {\s}{{\vec{s}}}
\newcommand {\qs}{\q-\s_{l}}

\newcommand {\Rvir}{R_{\rm vir}}       
\newcommand {\rmin}{r_{\rm min}}

\newcommand {\Muh}{\hat{M}}
\newcommand {\Pin}{P_{\infty}}
\newcommand {\Pinp}{P_{\infty,+}}
\newcommand {\Zi}{Z_{\infty}}
\newcommand {\Hh}{\hat{H}}
\newcommand {\Hi}{H_{\infty}}
\newcommand {\Hhi}{\hat{H}_{\infty}}
\newcommand {\Pt}{\Phi^{t}}
\newcommand {\Pit}{\Phi_{\infty}^{t}}
\newcommand {\Opm}{\Omega^{\pm}}
\newcommand {\Op}{\Omega^{+}}
\newcommand {\Om}{\Omega^{-}}
\newcommand {\B}[1]{\l\|#1\ri\|}
\newcommand {\vep}{\varepsilon} 
\newcommand {\LA}{\left\langle}
\newcommand {\RA}{\right\rangle}
\newcommand {\pa}{\partial}
\newcommand {\Ule}{U_{l}^{\vep}}
\newcommand {\Uhle}{{\hat{U}_{l}^{\vep}}}
\newcommand {\Eth}{{E_{\rm th}}}
\newcommand {\uk}{{\underline k}} 
\newcommand {\Po}{{\cal P}_E}

\begin{document}
\title {On the Integrability of the $n$--Centre Problem}
\author{Andreas Knauf\thanks{Mathematisches Institut,
Universit\"{a}t Erlangen-N\"{u}rnberg,
Bismarckstr.\ $1 \eh$, D--91054 Erlangen, Germany.
e-mail: knauf@mi.uni-erlangen.de}\and
Iskander A.\ Taimanov\thanks{Institute of Mathematics,
630090 Novosibirsk, Russia.
e-mail: taimanov@math.nsc.ru}}
\date{}
\maketitle
\begin{abstract}
It is known that for $n\geq3$ centres and positive energies the
$n$-centre problem of
celestial mechanics leads to a flow with a strange repellor and
positive topological entropy.

Here we consider the energies above some threshold and show:
Whereas for arbitrary $g>1$ independent integrals of Gevrey class $g$ 
exist,
no real-analytic (that is, Gevrey class 1) independent integral exists.

\end{abstract}
%
%
%
\section{Introduction}
%
In \cite{BT} the existence of a smoothly integrable geodesic flow on a
compact mani\-fold with positive topological entropy was established.
Positivity of topological entropy is seen as an indication of complex
dynamics, whereas integrability of a Hamiltonian flow is a metaphor for
its simplicity. So coexistence of these two properties may not have 
been
expected.
In fact, as we show here, such a coexistence takes place in natural
physical problems.

To be more specific,
in this note we consider the $n$--centre problem of celestial mechanics
in $d=2$ and $3$ dimensions.
We denote by $\s_k\in\bR^d,\ Z_k\in\bR\setminus {\{0\}}$ the location
resp.\
strength of the $k$--th centre, assuming $\s_k\neq \s_l$ for $k\neq l$.
Then the Hamiltonian function
\[\hH:T^*\hM\to\bR \qmbox{,} \hH \pq =\eh {\p\,}^2+V(\q),\]
with potential
\beq
V:\hM\to\bR\qmbox{,} V(\q) = -\sum_{k=1}^n\frac{Z_k}{\|\q-\s_k\|},
\Leq{pot}
on the cotangent bundle $T^*\hM$ of configuration space
\[\hM := \bR^d\setminus {\{\s_1,\ldots,\s_n\}}\]
generates a -- in general incomplete -- flow.

Denoting by $\hat{\omega}:=\sum_{i=1}^ddq_i\wedge dp_i\rstr_{T^*\hM}$
the
restricted canonical symplectic form, there exists a unique smooth
extension
\beq(P,\omega,H) \qmbox{of the hamiltonian
system}(T^*\hM,\hat{\omega},\hH)
\Leq{ext}
such that the flow $\Phi:\bR\times P\to P$ of $H$ is complete
(see \cite{Kn}, Thm. 5.1).

Concerning integrability of the flow, the following is known:

\begin{itemize}

\item
For $n=1$ this Hamiltonian system is integrable, with angular momentum
\[{\vec L}:P\to\bR^3 \qmbox{,} {\vec L}(x) =\l\{\begin{array}{ccc}
(\q-\s_1)\times \p&,&x= \pq,\ \q\neq \s_1\\
0&,&{\rm otherwise} \end{array}\ri.\]
for dimension $d=3$ being a real analytic constant of motion. For
$Z_1>0$ this is called the Kepler problem.

\item
For $n=2$ centres one introduces elliptic prolate coordinates to
analytically integrate
the flow $\Phi$, see e.g.\ \cite{Ar} or \cite{Th}.

\item
For $n\geq3$ centres and $d=2$ Bolotin showed in \cite{Bo} the
nonexistence of an analytic integral of the motion which 
is non--constant on an energy
shell $H^{-1}(E),\ E>0$, see also the discussion in Fomenko \cite{Fo}.

\item
For $d=3$ and a collinear configuration of centres the angular momentum
w.r.t.\ that axis is an additional constant of the motion, independent
of the number $n$ of centres. However, for $d=3$
it was proved by the first author \cite{Kn} for
sufficiently large energies $E > \Eth$
and by Bolotin and Negrini \cite{BN,BN2}
for nonnegative energies $E \geq 0$ that
the topological entropy of the flow,
restricted to the set $b_E$ of bounded orbits on $H^{-1}(E)$ is positive
(and
$h_{\rm top}(E)=0$ for $b_E=\emptyset$). Furthermore $h_{\rm top}(E)$ is zero
for
$n=1$ and $2$, and $h_{\rm top}(E)>0$ if $n\geq3$ and all centres being
attracting
or not more than two $\s_k$ being on a line (for collinear
configurations
with $Z_1,\ldots,Z_n<0$ one has $h_{\rm top}(E)=0$ for $E>0$).
\end{itemize}

In the present paper it is proved that

{\sl
\begin{itemize}

\item
for $d=2$, attracting centres ($Z_k >0$) and $E >0$,

\item
for $d=3$, arbitrary $Z_k \neq 0$, non-collinear configurations of
centres and $E > \Eth$ where the threshold energy level $\Eth$
is determined by the data $Z_k$, $\vec{s}_k$, $k=1,\dots,n$,

\end{itemize}

\noindent
the $n$-centre problem restricted onto the energy level $E$ admits
$d-1$ independent integrals of motion which are smooth (and moreover 
are
of the Gevrey class $g$ for any $g>1$) (see {\bf Thm.\ 
\ref{maintheorem}}).
Therefore the restricted problem is smoothly integrable.

On the other hand if the affine span of the (non-collinear)
centres equals $\bR^3$
then the restricted problem is not integrable in the real-analytic 
sense, 
that is Gevrey class 1 (see {\bf Thm.\ \ref{thm:ana}}).}

The article is based on the analysis of the
$n$-centre problem given
in the paper by the first author \cite{Kn} where, in particular, it is
shown that for high energies there are maps which relate scattering 
orbits
to their asymptotics given by scattering orbits of the Kepler
problem which is integrable (the necessary facts extracted
from \cite{Kn} are exposed in sections 2--3 of the present paper).
Such a relation of the $n$-centre problem with an integrable problem
leads to the
smooth integrability of the $n$-centre problem which is proved here.
Here we use a trick for constructing smooth first integrals from
discontinuous preserved quantities similar to others used in 
\cite{Bu,BT}

We would like to notice that in \cite{Bo}
the nonexistence of an additional analytic integral of motion
for the two-dimensional problem was derived by
Bolotin from Kozlov's theorem \cite{Ko} which reads that the geodesic 
flow of
(real)-analytic metric on a closed oriented surface of genus $g >1$ 
does not
admit an additional analytic first integral.

Although one can expect that the integrability of
a problem by analytic integrals of motion implies vanishing of the
topological entropy it is not proved until now. Therefore the results 
from
\cite{BN,Kn} do not imply the nonexistence of a complete family of
analytic first integrals. We prove that by an analysis of the set 
formed
by bounded orbits which supports the restricted flow with positive
topological entropy.

\section{Known Smoothness Results}
%
%
There are three basic types of motion: bounded, scattering and trapped,
corresponding to the disjoint subsets $b,s,t\subset P$ with
\[b^{\pm} := \{x\in P\mid \q\,(\pm\bR^+,x)\mbox{ is bounded
}\}\qmbox{,}
b:=b^+\cap b^-,\]
\[s^{\pm} := \{ x\in P \mid x\not\in b^{\pm}\mbox{ and }
H(x)>0\}\qmbox{,}
s:=s^+\cap s^-\]
and $t:=s^+\triangle s^-$. The orbits in $s$ go to spatial infinity
in both time directions, but because of the long range character of the
effective potential of strength $\Zi:=\sum_{k=1}^n Z_k$, one describes
these limits by comparison with (regularized) Kepler flow $\Pit:\Pin\ar
\Pin$
generated by the extension of
\beq
 \Hhi:T^*(\bR^d\setminus{\{0\}})\ar\bR\qmbox{,}
 \Hhi\pq := \eh \p^{\,2} - \frac{\Zi}{\B{\q}}
\Leq{free:hamiltonian}
Identifying the two phase spaces $P$ and  $\Pin$ outside a region
projecting
to a ball in configuration space which contains all singularities,
the {\em M\o ller transformations}
\beq
\Opm:\Pinp := \{ x\in \Pin\mid \Hi(x)>0\}\ar s^\pm
\mbox{ , }
\Opm := \lim_{t\ar\pm\infty} \Phi^{-t}\circ\Id\circ\Pit
\Leq{Moeller}
exist as pointwise limits, and are measure-preserving diffeomorphisms
(Thm.\ 6,3 and 6.5 of \cite{Kn}).
Similarly the asymptotic limits of the momentum
\beq
\p^{\pm}:s^{\pm}\ar\bR^{d}\qmbox{,}
\p^{\pm}(x_0) := \lim_{t\ar\pm\infty}\p\circ\Phi^t(x_0)
\Leq{asympto}
are smooth. Finally we define {\em time delay}
$\tau:s\ar\bR$ of a scattering state $x\in s$
by
\begin{eqnarray}\label{(6)}
\tau(x) := \lim_{R\ar\infty} & &\!\!\!\!\!\!\!\!
\int_{\bR} \Big(\sigma(R)\circ\Pt(x) -
\label{eq:TD:def}\\
& & \eh ( \sigma_{\infty}(R)\circ\Pit\circ\Op_{*}(x) +
              \sigma_{\infty}(R)\circ\Pit\circ\Om_{*}(x) )\Big) dt,
\nonumber
\end{eqnarray}
where $\sigma(R):P\ar \{0,1\}$ and $\sigma_{\infty}(R):\Pin\ar \{0,1\}$
are the characteristic functions $\sigma(R)(x) :=
\theta(R-\B{\q(x)})$ and similarly for $\sigma_{\infty}(R)$.
That asymptotic difference between the time  spent by the orbit and its
Kepler
limits inside a ball of large radius diverges near $b\cup t$.
However $\tau$ is  smooth, as the M\o ller transformations are.
%
\section{Analyticity Properties}
%
\begin{proposition}\label{prop1}
For a potential $V$ of the form (\ref{pot}) the following maps are
real-analytic:
\begin{enumerate}
\item
the flow $\Phi:\bR\times P\ar P$
\item
the M\o ller transformations
$\Opm:\Pinp\ar s^\pm$ and asymptotic momenta $\p^\pm: s^\pm\ar\bR^d$
\item
the time delay $\tau: s\ar\bR$.
\end{enumerate}
\end{proposition}
{\bf Proof.}\\
{\bf 1)}
We first indicate the definition of phase space $P$ in order to show
that for the potential (\ref{pot}) the smooth extension
(\ref{ext}) actually works in the real-analytic category
(in \cite{Kn} more general non-analytic potentials $V$ were
considered).
We assume $d=3$, the case of $d=2$ dimensions following by restriction.

For small $\vep>0$ in the phase space neighbourhood
\beq
\Uhle := \l\{\pq\in T^{*}\Muh \l| \B{\q-\s_{l}}<\vep,\,\,
|\p|^{2}>\frac{3}{2}\frac{Z_{l}}{\B{\q-\s_{l}}} \ri.\ri\}.
\Leq{U:l}
of the $l$th centre the following real-analytic coordinates are used:
\begin{itemize}
\item
The restriction of the Hamiltonian function,
which splits into
\beq
\Hh\pq = \Hh_{l}\pq + W_{l}(\q)\qmbox{with}
\Hh_{l}\pq := \eh \p^{\,2} - \frac{Z_{l}}{\B{\q-\s_{l}}}
\Leq{Zerlegung:Ham}
and
\[W_{l}(\q) = \sum_{i\neq l} \frac{-Z_{i}}{\B{\q-\s_{i}}}.\]
\item
The angular momentum
\beq
\hat{L}_{l}:\Uhle\ar \bR^3\qmbox{,}\hat{L}_{l}\pq := (\qs)\times\p .
\Leq{def:hat:Ll}
relative to the position $\s_l$.
\item
The time
$\hat{T}_{l}:\Uhle\ar\bR$
after which the Kepler orbit generated by $\Hh_{l}$ is in its
pericentre
w.r.t.\ $\s_l$
\beq
\hat{T}_{l}\pq := \int^{\B{\qs}}_{\rmin\pq}
\frac{r\,dr}{\sqrt{2r^{2}\Hh_{l}\pq + 2Z_{l}r -
\hat{L}_{l}^{2}\pq}}\cdot
\sign((\q-\s_{l})\cdot\p).
\Leq{peric:time}
The neighbourhood $\Uhle\subset T^{*}\Muh$ is defined in a way which
makes
the pericentre unique. Its distance from $\s_l$ equals
\beq
\rmin\pq :=
\l\{ \begin{array}{cc}

\frac{-Z_{l}+\sqrt{Z_{l}^{2}+2\Hh_{l}\pq\hat{L}_{l}^{2}\pq}}{2\Hh_{l}\pq}\!\!
         &, \Hh_{l} \neq 0 \\
     \hat{L}_{l}^{2}\pq/2Z_{l} \!\! & , \Hh_{l} = 0
     \end{array} \ri. .
\Leq{rmin}
\item
The Runge-Lenz vector relative to $\s_l$ is given by
\[
\vec{F}_{l}: \Uhle\ar \bR^3\qmbox{,}\vec{F}_{l}\pq :=
\p\times\hat{L}_{l}\pq - Z_{l}\frac{\qs}{\B{\qs}}.
\]
Since $|\vec{F}_{l}|^{2} >Z_{l}^{2}/4>0$, the {\em pericentral
direction}
\beq
\hat{F}_l:\Uhle\ar S^2\qmbox{,}\hat{F}_l:=\vec{F}_{l}/|\vec{F}_{l}|
\Leq{peri}
is well-defined.
\end{itemize}
As $\hat{L}_l\cdot \hat{F}_l = 0$, we get a real-analytic
diffeomorphism onto the image
\[\hat{{\cal Y}}:\Uhle\ar T^*(\bR\times S^2)\setminus \bar{0},\quad
\pq\mapsto (\hat{T}_{l},\hat{L}_l; \Hh_{l},\hat{F}_l)\]
onto a punctured neighbourhood of the zero section $\bar{0}$ of
the cotangent bundle.
The flow generated by $\Hh_{l}$ is linearized in these coordinates, and
extended on
$\Ule := \Uhle\cup \l(\bR \times S^{2}\ri)$ by
\[ {\cal Y} := (T_{l},{\bf L}_{l};H_{l},{\bf F}_{l}):
\Ule \ar T^*(\bR\times S^{2})
\mbox{ , } {\cal Y}(x):=\l\{
\begin{array}{lc}
(0,0;h,{\bf f})&(h,{\bf f}) \in \bR\times S^{2}\\
\hat{{\cal Y}}(x) & \mbox{ otherwise}
\end{array}\ri.\]
to a full neighbourhood of the zero section.

In a final step the hamiltonian flow $\Phi:\bR\times P\ar P$ generated
the  continuous extension $H:P\ar \bR$ of the Hamiltonian function
$\Hh: T^{*}\Muh\ar\bR$ is linearized by slightly changing the
coordinates
${\cal Y}$ using the (real-analytic but incomplete) flow generated by
$\Hh$.
See \cite{Kn} for details.

In summary, we obtain a real-analytic extension (\ref{ext}) of the
Hamiltonian
system $(T^*\hM,\hat{\omega},\hH)$. Thus the Hamiltonian vector field
$X_H$
(defined by $i_{X_H}\omega=dH$) is real-analytic, too, and thus
the flow
$\Phi:\bR\times P\ar P$ is known to be real-analytic (see, e.g.\
\cite{Ho}),
proving assertion 1.
\\[2mm]
{\bf 2)}
This, however is insufficient to show real-analyticity of
the M\o ller transformations and
asymptotic momenta. It is known that even for smooth potentials these
maps
may be very non-smooth, see \cite{Si}.

There exists an energy-dependent {\em virial radius} $\Rvir>0$, with
\beq
\frac{d}{dt} \LA\q,\p \RA > \frac{E}{2} > 0 \qquad\mbox{if }
 \B{\q}\geq \Rvir(E)\mbox{ and }E:=H\pq.
\Leq{qp}
If we assume
$q_0:=\|\vec{q}_0\|\geq\Rvir(E)\qmbox{and}\LA \q_0,\p_0\RA\geq 0$,
then
\beq
\|\q(t)\| \geq q_0\cdot \langle t\rangle_{\lambda}
\qmbox{for all $t\geq 0$, with}\lambda:= \sqrt{E/2}/q_0
\Leq{q:large}
and
\[\langle t\rangle_{\lambda}:=\sqrt{1+(\lambda t)^2}\qmbox{,}
\langle t\rangle:= \langle t\rangle_{1}.\]
In particular a trajectory leaving the ball of radius $\Rvir(E)$
cannot reenter this ball in the future but must go to spatial infinity.

Analyticity estimates for the trajectory are then derived from the
integral
equation with initial conditions $x_0=(\p_0,\q_0)$
\[\q(t,x_0)=\q_0+t\p_0 - \int_0^t\int_0^s \nabla V(\q(\tau,x_0))\,d\tau
\,ds,\]
using the decay property of the potential
\[\big|\pa^\beta V(\q)\big|\leq  |\beta|! \l(
\frac{C}{\|\q\|}\ri)^{|\beta|+1}
\qquad(\beta\in\bN_0^3)\]
outside the virial radius, and (\ref{q:large}).

Setting
\[\|\vec{w}\|_{\lambda}:=\sup_{t\geq0} \frac{\B{\vec{w}(t)}}{\langle
t\rangle_{\lambda}}\]
we obtain for $\gamma:=(\alpha,\beta)$,
$\pa^\gamma_{x_0}:=\pa^\alpha_{p_0}\pa^\beta_{q_0}$
\beq
\hspace*{-5mm}\|\pa^{\gamma}_{x_0}\q(\cdot,x_0)\|_{\lambda}\leq
\alpha!\, q_0^{-|\beta|+\delta_{|\gamma|,1}}
        E^{-\eh |\alpha|-1+\delta_{|\gamma|,1}}
.
\Leq{ind:ass}
This allows us to perform the (locally uniform w.r.t.\ $x_0$) time
limit
in
\beqn\label{(16)}
\lefteqn{\pa^\gamma_{x_0}(\p(t,x_0)-\p_0) = -}&&\label{oneminp}\\
&&\hspace*{-15mm}\sum_{N=1}^g
\sum_{\stackrel{\gamma^{(1)}+\ldots+\gamma^{(N)}=\gamma}{|\gamma^{(i)}|>0}}
\int_0^t
D^N \nabla V(\q(\tau,x_0))
\l(\pa^{\gamma^{(1)}}_{x_0}\q(\tau,x_0),\ldots,
\pa^{\gamma^{(N)}}_{x_0}\q(\tau,x_0)\ri)d\tau. \NN
\eeqn
with $g:=|\gamma|\geq 1$, and to conclude that $\p^{\,+}$ is
real-analytic
at $x_0$. We can substitute the assumption
$x_0\in s^+$ for the stronger assumptions
$q_0:=\|q_0\|\geq\Rvir(E)$, $\LA \q_0,\p_0\RA\geq 0$,
as initial conditions meeting the first lead to data meeting the
second one after some $t$.
The same holds for $\p^{\,-}$ by reversibility
($\p^{\,-}(\p_0,\q_0)= -\p^{\,+}(-\p_0,\q_0)$).

The proof of real-analyticity of the M\o ller transforms is
based on the integral equation
\[\pa^\gamma_{x_0}\vec{r}(t) = \int_{t}^{\infty} \int_{s}^{\infty}
\ \pa^\gamma_{x_0}\nabla
\l(\frac{-\Zi}{\B{\Q(\tau;x_0)}} - V(\q(\tau,x_0)) \ri)    d\tau\,  ds
\NN\\
\]
for $\vec{r}(t;x_0):=\q(t,x_0)-\Q(t;x_0)$ with Kepler hyperbola
$\l(\P(t;x_0),\Q(t;x_0)\ri) = \Pit(X_0)$.
Inspecting the proof of smoothness for the M\o ller transforms in
\cite{Kn},
one sees that the estimates for $\pa^\gamma_{x_0}\vec{r}$ can be
dominated by
$|\gamma|! C^{|\gamma|}$.\\[2mm]
{\bf 3)}
For a scattering state
$x\in s$ time delay equals
\beqno
\!\!\!\!\tau(x)
&=&\eh{\displaystyle \lim_{R\ar\infty}} \int_{\bR}
\sigma^{+}\circ\Pt(x)\cdot
\l(\sigma_{\infty}(R)\circ\Op_{*} -
\sigma_{\infty}(R)\circ\Om_{*}\ri) (\Pt(x))dt\\
&+&
\eh{\displaystyle \lim_{R\ar\infty}} \int_{\bR}
\sigma^{-}\circ\Pt(x)\cdot
\l(\sigma_{\infty}(R)\circ\Om_{*} -
\sigma_{\infty}(R)\circ\Op_{*}\ri) (\Pt(x))dt,
\eeqno
with $\sigma^{\pm}\pq := \theta(\pm\q\cdot\p)$, see \cite{KK} and
$\Opm_{*} :=
(\Opm)^{-1}$.
Using the intertwining property $\Opm_{*}\circ\Pt=\Pit\circ \Opm_{*}$
this can be written entirely in terms of M\o ller transformations and
the Kepler flow:
\beqno
\hspace*{-25mm}\tau(x)
&\hspace*{-3mm}=&\hspace*{-3mm}\eh{\displaystyle \lim_{R\ar\infty}}
\int_{\bR} \sigma^{+}\circ\Pit(x)\cdot
\l(\sigma_{\infty}(R)\circ\Pit\circ\Op_{*}(x) -
\sigma_{\infty}(R)\circ\Pit\circ\Om_{*}(x)\ri)dt \\
&\hspace*{-3mm}-&\hspace*{-3mm}
\eh{\displaystyle \lim_{R\ar\infty}} \int_{\bR}
\sigma^{-}\circ\Pit(x)\cdot
\l(\sigma_{\infty}(R)\circ\Pit\circ\Op_{*}(x) -
\sigma_{\infty}(R)\circ\Pit\circ\Om_{*}(x)\ri)dt .
\eeqno
With Assertion 2\ and the analog of formula (\ref{peric:time})
this implies analyticity of $\tau$. \hfill $\Box$\\
%
\section{Gevrey Integrals of Motion}
%
We now show the existence of independent constants of motion for all
energies
$E>\Eth$.
In \cite{Kn} many estimates are shown to hold true above a threshold
energy
$\Eth\geq0$:

In \cite{Kn} many estimates are shown to hold true above a threshold
energy
$\Eth\geq0$:

\begin{itemize}
\item
For $d=2$ and attracting centres ($Z_k>0$) $\Eth=0$.
\item
For $d=3$, arbitrary $Z_k\neq0$ and non-collinear configurations of the
$\vec{s}_k\in\bR^3$ the existence of such a threshold is proven.
\end{itemize}
In particular the set $b_E$ of bounded orbits of energy $E>\Eth$ is
shown
to be of measure zero (and has a Cantor set structure for $n\geq 3$).

According to the standard definition (see, e.g.\ \cite{AM}, Def.\
5.2.20)
functions $f_0,\ldots,f_k:\tilde{P}\to\bR$ on a symplectic manifold
$(\tilde{P},\tilde{\omega})$ are called {\em independent} if the set of
singular points of $F := f_0\times\ldots\times f_k:\tilde{P}\to\bR^k$
has
measure zero.

To simplify discussion, we set
\[\tilde{P} := H^{-1}([E_1,E_2])\]
for an
arbitrary energy interval, $\Eth\leq E_1\leq E_2<\infty$, and
$f_0^g\equiv
f_0:=H|_{\tilde{P}}$. Then for a parameter $g>1$ we define for $k=1,2$
\beq
f_k^g:\tilde{P}\to\bR \qmbox{,}
f_k^g(x) := \l\{\begin{array}{ccc}
p_k^+(x)\exp\l(-e^{{C(g)}\LA\tau(x)\RA}\ri)&,&x\in s\\
0&,&x\not\in s, \end{array}\ri.
\Leq{f}
where $C(g) = \frac{C_2}{g-1}$ with $C_2$ to be defined
below. Putting things together we get a map $F^g:\tilde{P}\to\bR^3$.

The notation is chosen so that $F^g$ belongs to the Gevrey class of
index
$g$.

Now we collect some information about Gevrey functions, which were
introduced in \cite{Ge}.
\begin{definition}
For $g\geq1$ and an open set $\Omega\subset\bR^n$ a function $f\in
C^\infty
(\Omega)$ is called of {\em Gevrey class} $g$ if for every compact
$K\subset\Omega$ there exist $A_K,C_K>0$ with
\[\max_{x\in K}|\pa^\alpha f(x)|\leq A_K\, C_K^{|\alpha|}\, (\alpha!)^g
\qquad
(\alpha\in\bN_0^n).\]
Here $|\alpha| = \alpha_1+\ldots+\alpha_n$ and
$\alpha!=\alpha_1!\cdot\ldots\cdot\alpha_n!$ if
$\alpha=(\alpha_1,\ldots,\alpha_n)$.\\
The vector space of these functions is denoted by $G_g(\Omega)$.
\end{definition}
Then
$G_1(\Omega)$ is the space of real-analytic functions, and
$G_{g'}(\Omega)
\supset G_g(\Omega)$ for $g'>g$.
$G_g(\Omega)$ is stable w.r.t.\ partial derivatives, compositions, and
the
implicit function theorem holds within the class. For $g>1$ Gevrey
partitions of unity exist.
\begin{example}
For $g>1$
\[f^g(x) := \l\{\begin{array}{ccc}\exp(-x^{-1/(g-1)})&,&x>0\\
0&,&x\leq0,\end{array}\ri.\]
is a function in $G_g(\bR)$ but not in $G_{g'}(\bR)$
for any $g'<g$, see e.g.\ \cite{Ju}.
\end{example}
For a real-analytic manifold we can define the spaces of Gevrey
functions,
as the defining family of bounds is preserved by coordinate changes.
\begin{theorem}
\label{maintheorem}
For all $g>1$ the functions $f_0^g,f_1^g,f_2^g$ are independent, in
involution
and of Gevrey class $g$.
\end{theorem}
{\bf Proof.}\\
In Thm.\ 12.8 of \cite{Kn} it was shown that the set $b_E$ of bounded
orbits of energy $E$ is of Liouville measure zero for all $E>\Eth$. As
the
set $t_E := \{x\in H^{-1}(E)\mid x\in T\}$ of trapped orbits is always
of
Liouville measure zero, on $\tilde{P}$ the set $\tilde{P}\setminus s$ 
is
of
measure zero. On $\tilde{P}\cap s$ the functions $f_1^g$ and $f_2^g$
are
real-analytic, using Proposition \ref{prop1}, whereas they are zero on
$\tilde{P}\backslash s$.

We study their decay near $\tilde{P}\setminus s$ in order to prove that
they are in $G^g(\tilde{P})$.

All orbits in $\tilde{P}\setminus s$ enter the interaction zone, so we
need only orbits in $\tilde{P}\cap s$ entering the interaction zone.
W.l.o.g.\
we assume $\Rvir$ to be constant on $[E_1,E_2]$. Then the region
\[\tilde{I} := \{x\in\tilde{P}\mid\|q(x)\|\leq\Rvir\}\]
projecting to the interaction zone is compact. The restriction of the
real-analytic flow $\Phi : \bR\times P\to\bR$ to a domain of the form
$[-\vep,\vep]\times\tilde{I}$ thus has partial derivatives $|\pa^\alpha
\Phi(t,x)|\leq \tilde{C}_1\tilde{C}_2^{|\alpha|}  \alpha! $.
We conclude, using Corollary \ref{cor} that for {\em arbitrary}
$t\in\bR$
such that $\Phi(t,x)\in\tilde{I}$, too,
\beq
|\pa^\alpha\Phi(t,x)|\leq C_1\exp({C_2}|\alpha|\,\LA t\RA)\ \alpha!.
\Leq{flow:e}
Next we analyse time delay (\ref{(6)}) for orbits entering $\tilde{I}$.
By compactness of $\tilde{I}$ there exists a $\tilde{R}\geq\Rvir$ such
that all Kepler orbits $\{\Phi_\infty^t\circ\Omega_*^\pm(x)\mid
t\in\bR\}$
through points $x\in\tilde{I}$ enter a configuration space region of
radius
$\tilde{R}$. In (\ref{(6)}) we only consider radii $R\geq\tilde{R}$
and,
denoting by $\tau_R(x)$ the integral (\ref{(6)}), we split this
quantity in
\[\tau_R(x) = \tau_R^+(x)+\tau_R^-(x)+\tau^0(x).\]
Here
\begin{eqnarray}\label{X}
\tau_R^\pm(x) &:=&
\int_\bR\big(\sigma^+\circ\Phi_t(x)\cdot\sigma(R)(1-\sigma
(\Rvir))\circ\Phi^t(x)\nonumber\\
&&\qquad-\eh\sigma_\infty(R)\circ\Phi_\infty^t\circ\Omega_*^\pm(x)\big)\,dt
\end{eqnarray}
and
\beq
\tau^0(x) := \int_\bR\sigma(\Rvir)\circ\Phi^t(x)\,dt.
\Leq{Y}
Whereas the limits $\tau^\pm(x):=\lim_{R\ar\infty}\tau_R^\pm(x)$
of (\ref{X}) meet Gevrey class 1 estimates uniform on $\tilde{I}\cap
s$,
$\tau^0(x)$ defined in (\ref{Y}),
the time spent by the orbit in the interaction zone, is not uniformly
bounded
on $\tilde{I}\cap s$ for $n\geq2$ centres, see Figure 1.
\begin{figure}
\begin{center}
\begin{picture}(0,0)%
\includegraphics{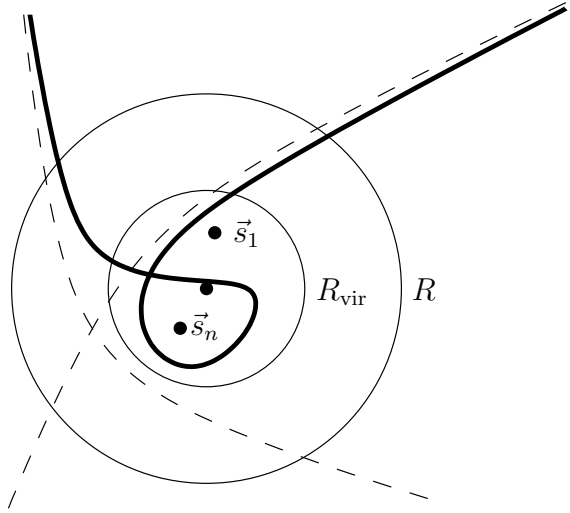}
\end{picture}%
\setlength{\unitlength}{3947sp}%
\setlength{\unitlength}{3947sp}%
\begingroup\makeatletter\ifx\SetFigFont\undefined%
\gdef\SetFigFont#1#2#3#4#5{%
  \reset@font\fontsize{#1}{#2pt}%
  \fontfamily{#3}\fontseries{#4}\fontshape{#5}%
  \selectfont}%
\fi\endgroup%
\begin{picture}(3562,3241)(859,-2428)
\put(2819,-1073){\makebox(0,0)[lb]{\smash{\SetFigFont{12}{14.4}{\rmdefault}{\mddefault}{\updefault}{$\Rvir$}%
}}}
\put(3420,-1073){\makebox(0,0)[lb]{\smash{\SetFigFont{12}{14.4}{\rmdefault}{\mddefault}{\updefault}{$R$}%
}}}
\put(2296,-706){\makebox(0,0)[lb]{\smash{\SetFigFont{12}{14.4}{\rmdefault}{\mddefault}{\updefault}{$\s_1$}%
}}}
\put(2026,-1291){\makebox(0,0)[lb]{\smash{\SetFigFont{12}{14.4}{\rmdefault}{\mddefault}{\updefault}{$\s_n$}%
}}}
\end{picture}
\end{center}
\caption{Scattering trajectory in configuration space (solid line),
with associated Kepler hyperbolae (broken lines)}

\end{figure}

We need only consider scattering states in a slightly smaller domain
$\tilde{I}_\vep := \{x\in\tilde{P}\mid\|q(x)\|\leq \Rvir-\vep\}$. By
(\ref{qp}) the configuration space trajectories with initial conditions
in $\tilde{I}_\vep\cap s$ intersect the boundary of the interaction
zone
with angles uniformly bounded away from zero. Therefore (\ref{flow:e})
shows
that for $x\in\tilde{I}_\vep\cap s$
\[|\pa^\alpha\tau^0(x)|\leq C_1 \exp(C_2|\alpha|\LA \tau(x)\RA)\
\alpha!,\]
and $\tau = \tau^+ + \tau^- + \tau^0$
satisfies the same kind of estimate (with enlarged constants).

A similar kind of reasoning applies to the asymptotic momenta $\p^\pm$:
\[|\pa^\alpha p_k^+(x)|\leq C_1\exp(C_2|\alpha|\LA\tau(x)\RA)\
\alpha!\]
for $x\in\tilde{I}_\vep\cap s$.

>From this one can conclude that $f_1^g$
and $f_2^g$ are of Gevrey class $g$ for $C(g)$ large enough, see
Proposition \ref{propo:gevrey}.

The functions $\p^\pm$ and $\tau$ are $\Phi_t$--invariant.
Therefore $f_1^g$ and $f_2^g$ are $\Phi_t$--invariant, too. As $f_0^g =
H
\setminus {\tilde{P}}$ generates the flow $\Phi_t$ on $\tilde{P}$, the
Poisson brackets $\l\{f_0^g,f_1^g\ri\}$ and $\l\{f_0^g,f_2^g\ri\}$
vanish.

By definition $\p^\pm(x) = \lim_{t\to\pm \infty}p\circ\Phi_t(x)$, and
the decay
estimates (\ref{(16)}) imply that
\[\l\{\lim_{t\to\infty}p_1\circ\Phi_t(x),\lim_{t\to\infty}p_2\circ\Phi_t
(x)\ri\} = \lim_{t\to\infty}\{p_1,p_2\}\circ\Phi_t(x) = 0.\]
Thus the Poisson bracket $\l\{f_1^g,f_2^g\ri\}$ vanishes, too.

To see that the functions $f_k^g$ are independent, it suffices to show
that the (analytic) restriction
\[F^g|_{\tilde{P}\cap s}\]
has a measure zero set of singular points.
As the energy range $[E_1,E_2]$ under consideration consists of regular
values of $H$, we may consider instead the independence of $f_1^g$ and
$f_2^g$, restricted to $H^{-1}(E)\cap s$, for $E\in[E_1,E_2]$.

Independence of $f_1^g$ and $f_2^g$, restricted to an energy surface
$H^{-1}(E)$ with $E\in[E_1,E_2]$ follows from their real-analyticity on
$H^{-1}(E)\cap s$ and the independence of $p_1^+$ and $p_2^+$,
remarking
that for $\LA \q,\p\RA\to\infty$, $\p^+(\pq)\to \p$.
\hfill $\Box$
\begin{proposition}\label{map}
Consider real-analytic maps
\[f:\Omega_f\to\bR^d \qmbox{,} g:\Omega_g\to\Omega_f\]
with open sets
$\Omega_f,\Omega_d\subset\bR^d$, meeting
\[|\pa^\alpha f_i|\leq M_f\,C_f^{|\alpha|}\,\alpha! \qmbox{,}
|\pa^\alpha g_i|
\leq M_g\,C_g^{|\alpha|}\,\alpha! \qquad (i=1\, , \,
\ldots,d,\alpha\in\bN_0^d).\]
Then
\[|\pa^\alpha(f\circ g)_i|\leq M_f\,C_{fg}^{|\alpha|}\,\alpha!
\qmbox{with}
C_{fg} := C_g(1+C_fdM_g)d.\]
\end{proposition}
{\bf Proof.}
The chain rule has the form
\[\pa^\alpha f\circ g =
\sum_{k=1}^n\sum_{A_1,\ldots,A_k}D^kf(\pa^{\alpha_
{A_1}}g,\ldots,\pa^{\alpha_{A_k}}g)\]
with $n=|\alpha|,\ (A_1,\ldots,A_k)$ running through the
$k$--partitions of
$\{1,\ldots,n\}$ and $\alpha_{A_l}\in\bN_0^d$ is the multiindex of size
$|\alpha_{A_l}| = |A_l|$ corresponding to $A_l\subset\{1,\ldots,n\}$ in
a, say
lexicographic, ordering of the $n$ partial derivatives.
Therefore with $\lambda:=C_fdM_g$
\beqno
|\pa^\alpha(f\circ g)| &\leq&
\sum_{k=1}^nM_fC_f^kd^kk!\sum_{A_1,\ldots,A_k}
\prod_{l=1}^kM_gC_g^{|A_l|}\alpha_{A_l}!\\
&=& M_fC_g^n\sum_{k=1}^n(C_fdM_g)^kk!\sum_{A_1,\ldots,A_k}\prod_{l=1}^k
\alpha_{A_l}!\\
&\leq&
M_fC_g^n\sum_{k=1}^n\lambda^kk!\sum_{A_1,\ldots,A_k}\prod_{l=1}^k
|A_l|!\\
&=& M_fC_g^nn!\sum_{k=1}^n\lambda^k{\textstyle {n-1\choose k-1}}\\
&=& M_fC_g^n\lambda(1+\lambda)^{n-1}n!,
\eeqno
since $k!\sum_{A_1,\ldots,A_k}\prod_{l=1}^k|A_l|! = n!{n-1\choose
k-1}$.
As in $d$ dimensions
\[|\alpha|! \leq d^{|\alpha|} \, \alpha! \qquad(\alpha\in\bN_0^d),\]
the estimate follows.
\hfill $\Box$\\[2mm]
This estimate is iterated in order to estimate the flow for long times:
\begin{corollary}\label{cor}
Assume that for $k=1,\ldots,t$ the maps
\[f^{(k)}:\Omega^{(k)}\to f^{(k)}
\big(\Omega^{(k)} \big)
\subset\Omega^{(k+1)}\qmbox{on}\Omega^{(k)}\subset\bR^d\]
admit the uniform estimates
$\big|\pa^\alpha f^{(k)}_i\big| \leq M\tilde{C}^{|\alpha|}
\,\alpha!$ with $M\geq1$.

\noindent
Then their iterates $T^{(k)} := f^{(k)}\circ T^{(k-1)}$ with $T^{(0)}
:=
\Id_{\Omega_1}$ are estimated by
\beq
\big|\pa^\alpha T^{(k)}_i\big|\leq
M\exp(C\,|\alpha|\,k)\,\alpha! \qquad
(\alpha\in\bN_0^d,k=1,\,\ldots,d)
\Leq{est}
with $C := \ln(d)+Md\tilde{C}$.
\end{corollary}
{\bf Proof.}
(\ref{est}) holds for $k=1$, and is assumed to hold for some $k<t$.
Then
by Proposition \ref{map}
\beqno
|\pa^\alpha T^{(k+1)}_i| &\leq& M(\exp(Ck)(1+Md\tilde{C})d)^{|\alpha|}
\,\alpha!\\
&\leq& M\exp(C(k+1))^{|\alpha|}\,\alpha!,
\eeqno
proving (\ref{est}).
\hfill $\Box$
\begin{proposition}\label{propo:gevrey}
For the choice $C(g) := \frac{C_2}{g-1}$ in (\ref{f}) the functions
$f_k^g$
are of Gevrey class $g$.
\end{proposition}
{\bf Proof.}
As $|\pa^\alpha\LA\tau(x)\RA|\leq
C_1\exp(C_2|\alpha|\LA\tau(x)\RA)\,\alpha!$,
\beqno
\lefteqn{|\pa^\alpha\exp(C(g)\LA\tau(x)\RA)| \leq
\exp\big(C(g)\LA\tau(x)\RA \big)\sum_{k=1}^{|\alpha|}
C(g)^k\sum_{A_1,\ldots,A_k}\prod_{l=1}^{k}\pa^{\alpha_{A_l}}\LA
\tau\RA}\\
&\leq& \exp\big((C(g)+C_2n)\LA\tau(x)\RA\big)
\l(\sum_{k=1}^nC(g)^kC_1^k{\textstyle{n-1\choose k-1}}\ri)\,\alpha!\\
&=& \exp\big((C(g)+C_2n)\LA\tau(x)\RA\big)(1+C(g)C_1)^n\,\alpha!
\eeqno
with $n := |\alpha|$.

The $\LA\tau\RA$ dependent part of $|\pa^\alpha f_k^g|$ is of the form
\[\exp\l(-e^{C(g)\LA\tau\RA}\ri)e^{C_2n\LA\tau\RA}\leq x^x,\]
as it has its maximum for $\LA\tau\RA$ with
\[x := e^{C(g)\LA\tau\RA} = {\textstyle\frac{C_2n}{C(g)}}.\]
Thus the choice $C(g) := \frac{C_2}{g-1}$ leads to the proof.
\hfill $\Box$

\section{Nonexistence of Analytic Integrals of Motion}
\label{nonanalytic}

We start with the following simple observation for a hamiltonian
system $(P,\omega,H)$: 

\begin{proposition} \label{depend}
Let $\gamma\subset \Sigma_E:=H^{-1}(E)$ be a periodic orbit of the 
Hamiltonian
flow $\Phi_t$ generated by the Hamiltonian function $H$ which is 
isolated on the
energy surface $\Sigma_E$. 
Then there is no additional integral of motion which is functionally 
independent
of $H$ on $\gamma$.
\end{proposition}

\noindent
{\bf Proof}. Let us assume that
there is an additional integral of motion $J:P\ar\bR$.
By localizing $J$ around $\gamma$ if necessary, the Hamiltonian flow 
$\Psi_s$  
generated by $J$ exists for all times $s\in\bR$.
Since these integrals are in involution, their flows commute:
$$
\Phi_t \Psi_s = \Psi_s \Phi_t\qquad (t,s\in\bR).
$$
Let $T$ be the period of $\gamma$ and $x \in \gamma$.
We have
$$
\Phi_{T} \Psi_s (x) = \Psi_s \Phi_{T} (x) = \Psi_s (x) \qquad 
(s\in\bR).
$$
This implies that $\Psi_s$ maps (within a given energy surface 
$\Sigma_E$)
periodic orbits of the flow $\Phi_t$ into periodic orbits of
this flow. But $\gamma\subset\Sigma_E$ 
is an isolated periodic orbit and therefore
$\Psi_s(x)\in\gamma$ for all $s\in\bR$. This
implies that $H$ and $J$ are functionally dependent on $\gamma$.
\hfill $\Box$
\begin{remark}
Note that {\em single} isolated periodic orbits need {\em not} 
form an obstruction to 
the existence of an additional analytic integral of motion $J:P\ar\bR$, 
independent of $H$ in the sense defined in Sect.\ 4. 
If the  periodic orbit $\gamma$ is hyperbolic, then
$J$ must be constant on its stable and unstable manifold, but still on 
a
neighbourhood of $\gamma\subset \Sigma_E$ the singular set of 
$J|_{\Sigma_E}$
may only consist of $\gamma$ which is of measure zero.

Easy examples for this are given by motion in a smooth potential 
$V:\bR^2\ar\bR$ which are rotationally symmetric. Then angular momentum
is an independent analytic integral of the motion, 
but circular hyperbolic orbits may exist for some energies. 

In the present context the energy surface $\Sigma_E\ (E>\Eth)$ for the 
two--centre
problem contains exactly one bounded orbit $\gamma$. $\gamma$ is 
periodic and 
collides with the two centres. As commented in the Introduction, the 
two--centre problem is analytically integrable.
\end{remark}
Let us return to the $n$-centre problem.
By Thm. 12.8 from \cite{Kn} the set of bounded orbits $b_E$ on the
energy level $E>\Eth$ consists in hyperbolic trajectories, and for 
$n\geq3$ 
is locally homeomorphic to a product of a Cantor set and the interval.
Moreover its Liouville measure vanishes. We now show that $b_E$ {\em 
can} 
form an obstruction to the existence of independent analytic integrals:

\begin{theorem} \label{thm:ana}
If the affine span of the (non-collinear) centres $\s_1,\ldots,\s_n$ 
equals $\bR^3$, then for
$E > E^\prime \geq \Eth$ 
the $n$--centre problem does not admit a pair of independent
analytic integrals $I_1,I_2:\Sigma_E\to\bR$ of motion 
on the energy surface $\Sigma_E:=H^{-1}(E)$.
\end{theorem}
{\bf Proof.} 
Assume that $I_1,I_2:\Sigma_E\to\bR$ are analytic integrals of motion. 
Then on the five--dimensional manifold
$\Sigma_E$ we consider the {\em singular set}
\[S:=\l\{x\in\Sigma_E\mid\rank(DI(x))<2\ri\}\]
with $I:=(I_1,I_2)$. 

By Prop.~\ref{depend} the periodic orbits within $b_E$ belong to the 
(closed) set $S$ (as in Prop.~\ref{depend} we consider functions 
$J:P\ar\bR$, we use appropriate smooth extensions $J$ of 
$I_k:\Sigma_E\to\bR$).

In \cite{Kn} it was shown that for energies $E>\Eth$ the 
periodic orbits are dense within $b_E$. Thus $b_E\subset S$, too.

According to \cite{vdD} (see also \cite{Ta} where a sketch
for more general subanalytic sets appears) $S$ admits an analytic 
simplicial
decomposition which is locally finite and whose simplices are 
semianalytic.
As $b_E\subset S$ is compact, there exist (disjoint) simplices 
$\Delta_1,
\ldots,\Delta_m$ with
\beq
b_E\subset\bigcup_{i=1}^m\Delta_i.
\Leq{fc}
Now we assume $I_1,I_2$ to be indendent integrals,
contrary to the statement of the theorem. Then
$K:=\max_i\dim(\Delta_i) < \dim(\Sigma_E)=5$.
As the set $S$ is $\Phi_t$--invariant, we consider the transversal 
intersections
\beq
\tilde{\Delta}_i:=\Delta_i\cap\cH_E
\Leq{dtilde}
of the simplices with the Poincar\'{e}
surface $\cH_E\subset\Sigma_E$ (defined in Sect.\ 10 of \cite{Kn}).
We have $\tilde{K}:=\max_i\dim(\tilde{\Delta}_i)=K-1<\dim(\cH_E)=4$.
The intersection $\Lambda_E:=b_E\cap\cH_E$ of the bounded orbits with 
the 
Poincar\'{e} surface has the form
\[\Lambda_E=\Lambda_E^+\cap\Lambda_E^-,\]
$\Lambda_E^\pm\subset\cH_E$ being the (un--) stable manifolds, 
consisting of
two--dimensional leaves which intersect transversally.

It is known from Hasselblatt \cite{Ha} that the H\"older regularity of 
(the distributions of) $\Lambda_E^\pm$ can be controlled by the 
so-called
bunching constant. 
For the case of the $n$-centre problem Prop.\ 11.2 of \cite{Kn} 
controls the expansion resp.\ contraction rates of the Poincar\'{e}
map on $\Lambda_E^\pm$, which differ from a constant times $E$ by 
$\cO(E^0)$. 
Thus applying \cite{Ha} we have 
$C^{2-\vep}$ regularity of $\Lambda_E^\pm$ for all large energies 
$E>\Eth$.

Now as $\tilde{K}<4$, for each $\tilde{\Delta}_i$ in (\ref{dtilde}) 
at least one of the intersections
\[\Delta_i^\pm:=\tilde{\Delta}_i\cap\Lambda_E^\pm\] 
must be of dimension $\leq1$. 
By reversibility of the flow we assume w.l.o.g.\ that 
$\dim(\Delta_1^-)\leq1$
and derive a contradiction.

As $\Lambda_E$ is a Cantor set (see Thm. 12.8 of \cite{Kn}), we assume 
w.l.o.g.\ 
that $\Delta_1^-$ contains a sequence $(x^{(i)})_{i\in\bN}$ 
of points $x^{(i)}\in\Lambda_E$
converging to $x\in\Delta_1^-$. By going to a subsequence, if 
necessary, we 
assume that the unit vectors $v^{(i)}\in T_x\cH_E$ with 
$\exp_x(t^{(i)}v^{(i)})=x^{(i)}$ have
a limit $v:=\lim_{i\to\infty}v^{(i)}$ (which is in fact independent of 
the choice of
Riemannian metric used to define the exponential map).

In \cite{Kn} starting from the alphabet $\cS:=\{1,\ldots,n\}$ the
space of symbol sequences 
\[{\bf X}:= \l\{{\underline k}\in \cS^\bZ \mid \forall\ i\in\bZ:
k_{i+1}\neq k_i \ri\}\]
is equipped with the metric
\[d({\underline k},{\underline l}) := 
\sum_{i\in\bZ}2^{-|i|}\cdot (1-\delta_{k_i,l_i}), \qquad
({\underline k},{\underline l}\in {\bf X}).\] 
Then denoting the {\em shift} by
\[\sigma:{\bf X}\ar{\bf X} \qmbox{,} \sigma({\underline k})_i := 
k_{i+1} 
\quad(i\in\bZ),\]
there exists a H\"older homeomorphism
\beq
\cF_E: {\bf X}\ar \Lambda_E
\Leq{co}
conjugating $\sigma$ and the restriction of the 
Poincar\'e map to $\Lambda_E$, see Lemma 12.2 of \cite{Kn}.

The  Poincar\'{e} surface is the disjoint union
\[\cH_E =\bigcup_{\stackrel{k,l=1}{ k\neq l}}^{n} \cH_E^{k,l},\]
the $\cH_E^{k,l}$ being open regions in the intersection of a 
five-dimensional
affine space and $\Sigma_E$. 
Starting with
\[V_E(k_{0},k_{1}):= W_E(k_{0},k_{1}):= \cH^{k_0,k_1}_E,\]
in  \cite{Kn} 
for $(k_{-m},\ldots,k_0)$ admissible (that is $k_l\neq k_{l+1}$) 
the nested subsets
\beq
W_E(k_{-m},\ldots,k_0):= 
W_E(k_{-1},k_0)\cap\Po(W_E(k_{-m},\ldots,k_{-1})),
\Leq{WEM}
resp.\ for $(k_{0},\ldots,k_m)$ admissible
\beq
V_E(k_{0},\ldots,k_m):= 
V_E(k_{0},k_1)\cap\Po^{-1}(V_E(k_{1},\ldots,k_{m})),
\Leq{VEM}
were defined, using the Poincar\'{e} map $\Po$.

By going to subsequences we can assume that the points $x^{(i)}\in 
\Lambda_E^-$ correspond to 
symbol sequences ${\underline k}^{(i)}\in {\bf X}$ which are related to 
the
symbol sequence ${\underline k}:=\cF_E^{-1}(x)\in {\bf X}$ of $x$ by
\[{k}^{(i)}_j= {k}_j \qquad (i\in\bN; j > \chi(i))\]
but ${k}^{(i)}_{\chi(i)}\neq  {k}_{\chi(i)}$,  where $\chi(i)\ar 
-\infty$. 

By construction the vector $v$ is tangent to the one-dimensional
manifold $\Delta_1^-$ at $x$.

We now show the existence of a second sequence $(y^{(i)})_{i\in\bN}$ of 
points
$y^{(i)}\in\Lambda_E$ converging to $x$, but with the following 
property: There
exists an $\alpha\in(0,1)$ such that writing the points in the form
\[y^{(i)}= \exp_x(s^{(i)}w^{(i)})\]
with units vectors $w^{(i)}\in T_x\cH_E$,
\beq
d(w^{(i)},{\rm span}(v))\geq |s^{(i)}|^{\alpha+1}.
\Leq{di}

Namely as the number $n$ of centres is $\geq 4$  (which
follows from our assumption on the positions $\s_1,\ldots,\s_n$)
we find symbol sequences ${\underline l}^{(i)}\in {\bf X}$ with
\[{l}^{(i)}_j= {k}_j \qquad (i\in\bN; j< \chi(i))\]
but ${l}^{(i)}_{\chi(i)}\neq  {k}_{\chi(i)}$, and 
\beq
{\rm affine\ span} 
\l(\s_{{k}_{\chi(i)-1}},\s_{{k}_{\chi(i)}},
\s_{{k}^{(i)}_{\chi(i)}},\s_{{l}^{(i)}_{\chi(i)}}\ri) =\bR^3.
\Leq{as}
Then the $y^{(i)}:=\cF_E(\underline{l}^{(i)})$ converge to
\[\lim_{i\to\infty}y^{(i)}=\cF_E\l(\lim_{i\to\infty}\underline{l}^{(i)}\ri)
=\cF_E(\uk)=x.\]
Next we consider the geometric situation at the (early) time $\chi(i)$.
More precisely we set
\[\tilde{\uk}:=\sigma^{\chi(i)}(\uk) \qmbox{,} 
\tilde{\uk}^{(i)}:=\sigma
^{\chi(i)}(\uk^{(i)}) \qmbox{and} \tilde{\ul}^{(i)}:=\sigma^{\chi(i)}
(\ul^{(i)}).\]
Thus
\[(\tilde{k}_1,\tilde{k}_2)=(\tilde{k}_1^{(i)},\tilde{k}_2^{(i)})=
(\tilde{l}_1^{(i)},\tilde{l}_2^{(i)})=\l(k_{\chi(i)+1},k_{\chi(i)+2}\ri)\]
but $\tilde{k}_0\neq \tilde{k}_0^{(i)} \neq \tilde{l}_0^{(i)} \neq 
\tilde{k}_0$.

By the conjugacy property of (\ref{co}) the identities 
$\tilde{x} := \cF_E(\tilde{\uk}) = \cP_E^{\chi(i)}(x)$, 
$\tilde{x}^{(i)} := \cF_E\l(\tilde{\uk}^{(i)}\ri) = 
\cP_E^{\chi(i)}(x^{(i)})$ 
and
$\tilde{y} := \cF_E\l(\tilde{\ul}^{(i)}\ri) = 
\cP_E^{\chi(i)}\l(y^{(i)}\ri)$ 
are true. 

All of these points are contained in the local stable manifold
\[V_E(\tilde{k}_1,\tilde{k}_2,\tilde{k}_3,\ldots)\subset 
V_E(\tilde{k}_1,\tilde{k}_2,\tilde{k}_3)\subset 
\cH_E^{\tilde{k}_1,\tilde{k}_2}, \]
but at the same time in the following disjoint sets 
\[\tilde{x}\in W_E(\tilde{k}_0,\tilde{k}_1,\tilde{k}_2)\qmbox{,} 
\tilde{x}^{(i)}\in 
W_E(\tilde{k}_0^{(i)},\tilde{k}_1,\tilde{k}_2)\qmbox{,} 
\tilde{y}^{(i)}\in W_E(\tilde{l}_0^{(i)},\tilde{k}_1,\tilde{k}_2).\]

In $\cH_E^{\tilde{k}_1,\tilde{k}_2}\subset\cH_E$ the minimal angle
between vectors from (points in) 
\[W_E(\tilde{k}_0,\tilde{k}_1,\tilde{k}_2)\cap 
V_E(\tilde{k}_1,\tilde{k}_2,\tilde{k}_3) 
\qmbox{to} W_E(\tilde{k}_0^{(i)},\tilde{k}_1,\tilde{k}_2)\cap 
V_E(\tilde{k}_1,\tilde{k}_2,\tilde{k}_3)\]
and vectors from
\[W_E(\tilde{k}_0,\tilde{k}_1,\tilde{k}_2)\cap 
V_E(\tilde{k}_1,\tilde{k}_2,\tilde{k}_3)
\qmbox{to} W_E(\tilde{l}_0^{(i)},\tilde{k}_1,\tilde{k}_2)\cap 
V_E(\tilde{k}_1,\tilde{k}_2,\tilde{k}_3)\]
is bounded away from zero by some $\Delta>0$ (see (\ref{VEM}) and 
(\ref{WEM})). 
This follows from (\ref{as}).
As there are only finitely many
(at most $n$, to be more precise)
choices for each of the indices 
$\tilde{k}_0,\tilde{k}_1,\tilde{k}_2,\tilde{k}_0^{(i)}$ and
$\tilde{l}_0^{(i)}$, this bound is uniform.

So the angle
between the unit vectors $\tilde{v}^{(i)}$ and $\tilde{w}^{(i)}$ 
defined by
\[\tilde{x}^{(i)}=\exp_{\tilde{x}}(\tilde{s}^{(i)}\tilde{v}^{(i)})\qmbox{,} 
\tilde{y}^{(i)}=\exp_{\tilde{x}}(\tilde{t}^{(i)}\tilde{w}^{(i)})\]
is bounded below by $\Delta $. 
But by Prop.\ 11.2 of \cite{Kn} the contraction rates of these vectors
w.r.t.\ the Poincar\'{e} $\cP_E$ differ at most by the order 
$\cO(1/E)$.

By estimate (11.2) of \cite{Kn} we conclude that after iterating 
$-\chi(i)>0$
times these vectors still have an angle bounded away from zero by
\[\Delta \l(1-\frac{C}{E}\ri)^{|\chi(i)|}\]
whereas their length is reduced at least by a factor 
$(cE)^{|\chi(i)|}$.

Thus by enlarging the energy $E$,
we can choose an arbitrarily small 
$\alpha >0$ in (\ref{di}).

But (\ref{di}) is incompatible with $y^{(i)}$ belonging to the 
one--dimensional
submani\-fold $\Delta_1^-$.

As we can apply the same argument to all $x\in\Lambda_E$, we have 
derived a 
contradiction to the finite covering assumption (\ref{fc}).
\hfill $\Box$

The second author (I.A.T.) was supported by RFBR (grant 03-01-00403) and
Max-Planck-Institute on Mathematics in Bonn.

\end{document}